\documentclass[11pt,reqno]{amsart}

\let\savethebibliography=\thebibliography
\usepackage[square,numbers,comma,sort&compress]{natbib}
\let\thebibliography=\savethebibliography

\usepackage{rotating}
\usepackage{graphicx}
\usepackage{psfrag}
\usepackage{amsfonts}
\usepackage{url}
\usepackage{color}
\usepackage{amssymb}
\usepackage{enumerate}

\newtheorem{theorem}{Theorem}

\newtheorem{Lemma*}[theorem]{Lemma}

\newtheorem{corollary}[theorem]{Corollary}

\theoremstyle{definition}
\newtheorem{definition}[theorem]{Definition}
\newtheorem{remark}[theorem]{Remark}

\def\jdlqed{\vbox{\hrule \hbox{\vrule\hbox to
5pt{\vbox to 6pt{\vfil}\hfil}\vrule}\hrule}}

\newcommand{\Q}{{\mathbb Q}}
\newcommand{\Z}{{\mathbb Z}}
\newcommand{\N}{{\mathbb N}}

\newcommand{\R}{\mathbb R}
\DeclareMathOperator\interior{int}
\DeclareMathOperator\vmin{vmin}
\newcommand\FPTAS{{\small FPTAS}}
\newcommand\VPareto{{V^{\mathrm{Pareto}}}}
\DeclareMathOperator\cone{cone}
\DeclareMathOperator\lcm{lcm}

\usepackage{ifthen}
\makeatletter
\newcommand{\DeclareBracket}[3]{
  \newcommand{#1}[2][]{%
  \ifthenelse%
  {\equal{##1}{}}%
  {\left#2##2\right#3}%
  {\csname ##1l\endcsname#2##2\csname ##1r\endcsname#3}}}    
\makeatother
\DeclareBracket\abs||           
\DeclareBracket\norm\|\|        
\DeclareBracket\floor\lfloor\rfloor
\DeclareBracket\ceil\lceil\rceil
\DeclareBracket\set\{\}
\DeclareBracket\paren()
\DeclareBracket\bracket[]
\DeclareBracket\inner\langle\rangle
\DeclareBracket\fractional\{\}
\makeatother

\usepackage{amssymb}

\usepackage{enumerate}
\def\ve#1{\mathchoice{\mbox{\boldmath$\displaystyle\bf#1$}}
{\mbox{\boldmath$\textstyle\bf#1$}}
{\mbox{\boldmath$\scriptstyle\bf#1$}}
{\mbox{\boldmath$\scriptscriptstyle\bf#1$}}}

\usepackage{ifpdf}

\usepackage[breaklinks=true,colorlinks,citecolor=blue,linkcolor=blue]{hyperref}

\title[Multicriteria Integer Linear Programs]{Pareto Optima of\\ Multicriteria Integer Linear Programs}

\author[J.~A.~De Loera]{Jes\'us~A.~De Loera}
\address{J.~A.~De Loera:
  University of California, Davis, Department of Mathematics,
  One Shields Avenue, Davis CA 95616, USA}
\thanks{J.~A.~De Loera gratefully acknowledges support from NSF grant
  DMS-0608785.}
\email{deloera@math.ucdavis.edu}

\author[R.~Hemmecke]{Raymond Hemmecke}
\address{R.~Hemmecke: Otto-von-Guericke-Universit\"at Magdeburg,
  Department of
  Mathematics, Institute for Mathematical Optimization (IMO),
  Universit\"ats\-platz~2, 39106 Magdeburg, Germany}
\email{hemmecke@imo.math.uni-magdeburg.de}

\author[M.~K\"oppe]{Matthias K\"oppe}
\address{M.~K\"oppe: Otto-von-Guericke-Universit\"at Magdeburg, Department of
  Mathematics, Institute for Mathematical Optimization (IMO),
  Univer\-si\-t\"ats\-platz~2, 
  39106 Magdeburg, Germany} 
\email{mkoeppe@imo.math.uni-magdeburg.de}
\thanks{M.~K\"oppe was supported by a Feodor Lynen Research Fellowship from the
  Alexander von Humboldt Foundation.}

\date{$\relax$Revision: 1.29 $ - \ $Date: 2007/07/09 22:49:05 $ $}

\begin{document}

\begin{abstract}
  We settle the computational complexity of fundamental questions related to
  multicriteria integer linear programs, when the dimensions of the strategy
  space and of the outcome space are considered fixed constants.  In
  particular we construct:
  \begin{enumerate}[1.]
  \item polynomial-time algorithms to exactly determine the number of Pareto
    optima and Pareto strategies;
  \item a polynomial-space
    polynomial-delay prescribed-order enumeration algorithm for arbitrary
    projections of the Pareto set;
  \item an algorithm to minimize the distance of a Pareto
    optimum from a prescribed comparison point with respect to arbitrary
    polyhedral norms;
  \item a fully polynomial-time approximation scheme for the problem
    of minimizing the distance of a Pareto optimum from a prescribed
    comparison point with respect to the Euclidean norm.
  \end{enumerate}
\end{abstract}

\maketitle

\section{Introduction}


\label{intro} 
Let $A = (a_{ij})$ be an integral $m \times n$-matrix and $\ve b \in
\Z^m$ such that the convex polyhedron $ P =  \{\, \ve u \in
\R^n  : A \ve u \leq \ve b \,\}\,$ is bounded. Given $k$~linear functionals
$f_1,f_2,\dots,f_k \in \Z^n$,
we consider the \emph{multicriterion integer linear programming problem}
\begin{equation}
  \label{eq:multicrit-ilp}
\begin{aligned}
 \vmin\quad& \bigl(f_1(\ve u),f_2(\ve u),\dots,f_k(\ve u)\bigr)\\
\hbox{subject to}\quad &  A\ve u \leq \ve b \\
& \ve u \in \Z^n
\end{aligned}
\end{equation}
where $\vmin$ is defined as the problem of finding all Pareto optima and a
corresponding Pareto strategy. For a lattice point $\ve u$ the vector
$\ve f(\ve u) = \bigl(f_1(\ve u),\dots,f_k(\ve u)\bigr)$ is called an \emph{outcome
vector}.  Such an outcome vector 
is a \emph{Pareto optimum} for the above problem if and only if there is no
other point~$\ve{\tilde u}$ in the feasible set such that $f_i(\ve{\tilde
  u})\leq f_i(\ve u)$ for all~$i$ and $f_j(\ve{\tilde u})
< f_j(\ve u)$ for at least one index $j$.  The corresponding feasible
point~$\ve u$ is
called a \emph{Pareto strategy}.  Thus a feasible vector is a Pareto strategy
if no feasible vector can decrease some criterion without causing a
simultaneous increase in at least one other criterion. For general information
about the multicriteria problems see, e.g., \cite{Figueiraetal:multicriteria-book,Sawaragietal:multiobjectivebook}.

In general multiobjective problems the number of Pareto optimal
solutions may be infinite, but in our situation the number of Pareto
optima and strategies is finite. There are several well-known
techniques to generate Pareto optima. Some popular methods used to
solve such problems include, e.g., weighting the objectives or using
a so-called global criterion approach (see \cite{ehrgott-gandibleux-2000}). In
abnormally nice situations, such as multicriteria \emph{linear}
programs~\cite{isermann-1974}, one knows a way to generate all Pareto optima, 
but most techniques reach only some of the Pareto optima.

The purpose of this article is to study the sets of \emph{all} Pareto optima and
strategies of a multicriterion integer linear program using the
algebraic structures of generating functions. The set of Pareto points
can be described as the formal sum of monomials
\begin{equation} \label{theparetogf}
 \sum \bigl\{\, \ve z^{\ve v} : 
 \text{  $\ve u \in P \cap \Z^n$ 
   and $\ve v = \ve f(\ve u) \in \Z^k$ is a Pareto optimum} \,\bigr\}.
\end{equation}
Our main theoretical result states that, under the assumption that the
number of variables is fixed, we can compute in polynomial time a
compact expression for the huge polynomial above, thus 
\emph{all} its Pareto optima can in fact be counted exactly. The same can
be done for the corresponding Pareto strategies when written in the form
\begin{equation} \label{theparetostgf}
 \sum \bigl\{\, \ve x^{\ve u} : \text{ $\ve u \in P \cap \Z^n$ and $\ve f(\ve
   u)$ is a Pareto optimum} \, \bigr\}.
\end{equation}
\begin{theorem} \label{mainsec1} 
Let $A \in \Z^{m \times n}$, a $d$-vector $\ve b$, and linear functions
$f_1,\dots,f_k\in\Z^n$ be given. There are algorithms to perform the following tasks:

\begin{enumerate}[\rm (i)]
\item Compute the generating function \eqref{theparetogf} of all the
Pareto optima as a sum of rational functions. In particular we can count 
how many Pareto optima are there.
If we assume $k$ and $n$ are fixed, the algorithm runs in time
polynomial in the size of the input data. \label{mainsec1-part-1}

\item Compute the generating function \eqref{theparetostgf} of all the
Pareto strategies as a sum of rational functions. In particular we can
count how many Pareto strategies are there in $P$. If we assume $k$
and $n$ are fixed, the algorithm runs in time polynomial in the size
of the input data. \label{mainsec1-part-2}

\item Generate the full sequence of Pareto optima ordered
  lexicographically or by any other term ordering.  If we assume $k$
  and $n$ are fixed, the algorithm runs in polynomial time on the input
  size and the number of Pareto optima.  (More strongly, there exists
  a \emph{polynomial-space polynomial-delay prescribed-order
  enumeration algorithm}.) \label{mainsec1-part-3}
\end{enumerate}

\end{theorem}

In contrast it is known that for non-fixed dimension it is $\#\mathrm{P}$-hard
to enumerate Pareto optima and $\mathrm{NP}$-hard to find them
\cite{emelichev-perepelitsa-1992,sergienko-perepelitsa-1991}.  
The proof of \autoref{mainsec1} parts (\ref{mainsec1-part-1}) and (\ref{mainsec1-part-2}) will be given in \autoref{theory}. It
is based on the theory of rational generating functions 
as developed in \cite{bar,BarviPom}. Part (\ref{mainsec1-part-3}) of \autoref{mainsec1} will
be proved in \autoref{s:listing}.\smallbreak

For a user that knows some or all of the Pareto optima or strategies,
a goal is to select the ``best'' member of the family. 
One is interested in selecting one Pareto optimum that realizes the
``best'' compromise between the individual objective functions.
The quality of the compromise is often measured by the distance of a
Pareto optimum~$\ve v$ from a user-defined comparison point~$\ve{\hat v}$.
For example, often users take as a good comparison point 
the so-called \emph{ideal point} $\ve v^{\textrm{ideal}}\in\Z^k$ of the
multicriterion problem, which is defined as 
    \begin{displaymath}
      v^{\mathrm{ideal}}_i = \min \{\, f_i(\ve u) : \ve u\in P\cap\Z^n \,\}.
    \end{displaymath}
The criteria of comparison with the point $\ve{\hat v}$ are quite diverse, but
some popular ones include computing the minimum over the possible sums of
absolute differences of the individual objective functions, evaluated
at the different Pareto strategies, from the comparison
point~$\ve{\hat v}$, i.e., 
\begin{subequations}
  \label{eq:global-criteria}
\begin{equation}
  f(\ve u) = |f_1(\ve u)-\hat v_1| +\dots +
|f_k(\ve u)-\hat v_k|,
\end{equation}
or the maximum of the absolute differences, 
\begin{equation}
  f(\ve u) = \max\bigl\{|f_1(\ve u)-\hat v_1|, \dots, |f_k(\ve u)-\hat v_k| \bigr\},
\end{equation}
over all Pareto optima $(f_1(\ve u),\dots,f_k(\ve u))$.  Another
popular criterion, sometimes called the \emph{global criterion}, is to
minimize the sum of 
relative distances of the individual objectives from their known minimal
values, i.e., 
\begin{equation}
  f(\ve u) = \frac{f_1(\ve u)-v^{\mathrm{ideal}}_1}{|v^{\mathrm{ideal}}_1|} + \dots + \frac{f_k(\ve u)-v^{\mathrm{ideal}}_k}{|v^{\mathrm{ideal}}_k|}.
\end{equation}
\end{subequations}

We stress that if we take any one of these functions as an objective
function of an integer program, the optimal solution will be a
non-Pareto solution of the multicriterion
problem~\eqref{eq:multicrit-ilp} in general.  In contrast, we show in
this paper that by encoding Pareto optima and strategies 
as a rational function we avoid this problem, since we evaluate the
objective functions directly on the space of Pareto optima.

All of the above criteria~\eqref{eq:global-criteria} measure the distance from
a prescribed point with respect to a \emph{polyhedral norm}. 
In \autoref{s:selection}, we prove:
\begin{theorem}\label{th:intro-polyhedral-opt}
  Let the dimension $n$ and the number  $k$ of objective functions be fixed.
  Let a multicriterion integer linear 
  program~\eqref{eq:multicrit-ilp} be given.  Let a polyhedral norm~$\mathopen\| \cdot
  \mathclose\|_Q$ be given by the vertex or inequality description of its unit ball
  $Q\subseteq\R^k$. 
  Finally, let a prescribed point $\ve{\hat v}\in\Z^k$ be given.

  \begin{enumerate}[\rm (i)]
\item There exists a polynomial-time algorithm to find a Pareto optimum~$\ve v$
  of~\eqref{eq:multicrit-ilp} that  
  minimizes the distance $ \| \ve v - \ve{\hat v}
  \|_Q$ from the prescribed point.
\item There exists a polynomial-space polynomial-delay enumeration algorithm
  for enumerating the Pareto optima of~\eqref{eq:multicrit-ilp} in the order
  of increasing distances from the prescribed point~$\ve{\hat v}$.
\end{enumerate}
\end{theorem}

Often users are actually interested in finding a Pareto optimum that minimizes
the \emph{Euclidean} distance from a prescribed comparison point~$\ve{\hat v}$, 
\begin{equation}
   f(\ve u) = \sqrt{{|f_1(\ve u)-\hat v_1|}^2 +\dots +
     {|f_k(\ve u)-\hat v_k|}^2},
\end{equation}
but to our knowledge no method of the literature gives a satisfactory
solution to that problem.  In \autoref{s:selection}, however, we prove the 
following theorem, which gives a very strong approximation result.
\begin{theorem}\label{th:intro-euclidean-fptas}
  Let the dimension~$n$ and the number~$k$ of objective functions be fixed.
  There exists a fully polynomial-time approximation scheme for the problem of
  minimizing the Euclidean distance of a Pareto optimum
  of~\eqref{eq:multicrit-ilp} from a prescribed comparison point $\ve{\hat
    v}\in\Z^k$.
\end{theorem}
We actually prove this theorem in a somewhat more general setting, using an
arbitrary norm whose unit ball is representable by a homogeneous polynomial
inequality. 

\section{The rational function encoding of all Pareto optima}  \label{theory}

We give a very brief overview of the theory of rational generating
functions necessary to establish \autoref{mainsec1}. For full details
we recommend \cite{bar,BarviPom,barvinok-woods-2003,deloera-hemmecke-koeppe-weismantel:intpoly-fixeddim} and the references therein.
In 1994 Barvinok gave an algorithm
for counting the lattice points in $P=  \{\, \ve u \in
\R^n  : A \ve u \leq \ve b \,\}\,$ in polynomial time when the
dimension~$n$ of the feasible polyhedron is a constant \cite{bar}.  The input for
Barvinok's algorithm is the binary encoding of the integers $a_{ij}$
and $b_i$, and the output is a formula for the multivariate generating
function 
\begin{displaymath}
  g(P; \ve x)=\sum_{\ve u \in P \cap \Z^n} \ve x^{\ve u}
\end{displaymath}
where $\ve
x^{\ve u}$ is an
abbreviation of $x_1^{u_1} x_2^{u_2}\dots x_n^{u_n}$. This long
polynomial with exponentially many monomials is encoded as a much
shorter sum of rational functions of the form
\begin{equation}
\label{barvinokseries} g(P; \ve x) = \sum_{i \in I} \gamma_i
\frac{\ve x^{\ve c_i}}{(1-\ve x^{\ve d_{i1}})(1-\ve x^{\ve d_{i2}})\dots
(1-\ve x^{\ve d_{in}})}.
\end{equation}
Barvinok and Woods in 2003 further developed a set of powerful manipulation
rules for using these short rational functions in Boolean
constructions on various sets of lattice points. 

Throughout the paper we assume that the polyhedron $ P =  \{\, \ve u \in
\R^n  : A \ve u \leq \ve b \,\}\,$ is bounded. We begin by recalling some useful results of Barvinok and Woods
(2003):

\begin{theorem}[Intersection Lemma; Theorem 3.6 in \cite{barvinok-woods-2003}] \label{intersect}
Let $\ell$ be a fixed integer.
Let $S_1, S_2$ be finite subsets of $\Z^n$. Let
$g(S_1; \ve x)$ and $g(S_2; \ve x)$ be their
 generating functions, given as short
rational functions with at most $\ell$ binomials in each denominator.
Then there exists a polynomial time algorithm, which computes
$$ g(S_1 \cap S_2; \ve x)  = \sum_{i \in I} \gamma_i \frac { \ve x^{\ve c_i} }
{  (1-\ve x^{\ve d_{i1}})  \dots (1-\ve x^{\ve d_{is}}) }$$
with $s \leq 2\ell$, where the $\gamma_i$ are rational numbers,
$\ve c_i,\ve d_{ij}$ are nonzero integer vectors, and $I$ is a polynomial-size
index set.
\end{theorem}

The following theorem was proved by Barvinok and Woods using
\autoref{intersect}: 

\begin{theorem}[Boolean Operations Lemma; Corollary 3.7 in \cite{barvinok-woods-2003}] \label{unioncomplement}
Let $m$ and $\ell$ be fixed integers.
Let $S_1, S_2, \dots,S_m$ be finite subsets of $\Z^n$. 
Let $g(S_i; \ve x)$ for $i=1, \dots, m$ be their generating
functions, given as short rational functions with at most $\ell$
binomials in each denominator.  
Let a set $S\subseteq\Z^n$ be defined as a Boolean combination of
$S_1,\dots,S_m$ 
(i.e., using any of the operations $\cup$, $\cap$, $\setminus$).
Then there exists a polynomial
time algorithm, which computes
$$ g(S; \ve x) =  \sum_{i \in I} \gamma_i \frac { \ve x^{\ve c_i} } {  (1-\ve
  x^{\ve d_{i1}})  \dots (1-\ve x^{\ve d_{is}}) }$$
where $s = s(\ell,m)$ is a constant, the $\gamma_i$ are rational numbers,
$\ve c_i,\ve d_{ij}$ are nonzero integer vectors, and $I$ is a polynomial-size index
set. 
\end{theorem}

We will use the \emph{Intersection Lemma} and the \emph{Boolean
Operations Lemma} to extract special monomials present in the
expansion of a generating function. The essential step in the
intersection algorithm is the use of the {\em Hadamard product}
\cite[Definition 3.2]{barvinok-woods-2003} and a special monomial substitution.
The Hadamard product is a bilinear operation on rational functions
(we denote it by $*$). The computation is carried out for pairs of
summands as in~\eqref{barvinokseries}.  Note that the Hadamard
product $m_1 * m_2$ of two monomials $m_1,m_2$ is zero unless
$m_1=m_2$. 

Another key subroutine introduced by Barvinok and Woods is the
following \emph {Projection Theorem}. 

\begin{theorem}[Projection Theorem; Theorem 1.7 in \cite{barvinok-woods-2003}] \label{project}
Assume the dimension $n$ is a fixed constant. Consider  a rational
polytope $P \subset \R^n$ and a linear map $T\colon \Z^n \rightarrow
\Z^k$. There is a polynomial time algorithm which computes a short
representation of the generating function $f \bigl(T(P \cap
\Z^n);\ve x\bigr) $.
\end{theorem}

One has to be careful when using earlier Lemmas (especially the Projection
Theorem) that the sets in question are finite. 
The proof of \autoref{mainsec1} will require us to  project
and intersect sets of lattice points represented by rational
functions. We cannot, in principle, do those operations for
\emph{infinite} sets of lattice points. Fortunately, in our setting it
is possible to restrict our attention to finite sets. 

Finally, one important comment. If we want to count the points of a
lattice point set~$S$, such as the set of Pareto optima, it would apparently
suffice to substitute 
$1$ for all the variables~$x_i$ of the generating function 
\begin{displaymath}
  g(S; \ve x) = \sum_{\ve u \in S} \ve x^{\ve u} = \sum_{i \in I} \gamma_i
\frac{\ve x^{\ve c_i}}{(1-\ve x^{\ve d_{i1}})(1-\ve x^{\ve d_{i2}})\dots
(1-\ve x^{\ve d_{in}})}
\end{displaymath}
to get the specialization $|S| = g(S;\ve x\,{=}\,\ve1)$. 
But this cannot 
be done directly due to the singularities in the rational function 
representation. Instead,
choose a generic vector $\ve \lambda=(\lambda_1,\dots,\lambda_n)$ and
substitute each of the
variables~$x_i$ by~$\mathrm e^{t\lambda_i}$. Then we get
\begin{displaymath}
g(S,\ve e^{t\ve \lambda})
= \sum_{i \in I} \gamma_i
\frac{\mathrm e^{t \langle \ve\lambda, \ve c_i\rangle}}
{(1-\mathrm e^{t \langle \ve\lambda, \ve d_{i1}\rangle})
  (1-\mathrm e^{t \langle \ve\lambda, \ve d_{i2}\rangle})\dots
  (1-\mathrm e^{t \langle \ve\lambda, \ve d_{in}\rangle})}.
\end{displaymath}
Counting the number of lattice points is
the same as computing the constant terms of the Laurent series for each
summand and adding them up. This can be done using elementary 
complex residue techniques (see \cite{BarviPom}).

\begin{proof}[Proof of \autoref{mainsec1}, part (\ref{mainsec1-part-1}) and (\ref{mainsec1-part-2})]
The proof of part~(\ref{mainsec1-part-1}) has three steps:\par\smallskip 
\noindent\emph{Step 1.} 
For $i=1,\dots,k$ let
$\bar v_i\in\Z$ be an upper bound of polynomial encoding size for the value of
$f_i$ over~$P$.  Such a bound exists because of the boundedness of~$P$, and it
can be computed in polynomial time by linear programming.  
We will denote the vector of
upper bounds by~$\ve{\bar v}\in\Z^k$. 
We consider the \emph{truncated multi-epigraph} of the objective functions
$f_1,\dots,f_k$ over 
the linear relaxation of the feasible region~$P$,
\begin{equation}
  \label{eq:multi-epigraph}
  \begin{aligned}
    P^\geq_{f_1,\dots,f_k} = \bigl\{\, (\ve u, \ve v)\in\R^n\times\R^k :{}
    &\ve u\in P,\\
    &\bar v_i \geq v_i \geq f_i(\ve u)\text{ for $i=1,\dots,k$} \,
    \bigr\},
  \end{aligned}
\end{equation}
which is a rational convex polytope in $\R^n\times\R^k$.  Let $V^\geq\subseteq\Z^k$
denote the integer projection of $P^\geq_{f_1,\dots,f_k}$ on the $\ve v$
variables, i.e., the set
\begin{equation}
  \label{eq:upper-value-set}
  V^\geq = \bigl\{\, \ve v \in \Z^k :
  \exists \ve u\in \Z^n \text{ with } 
  (\ve u,\ve v) \in P^\geq_{f_1,\dots,f_k} \cap(\Z^n\times\Z^k) 
  \,\bigr\}.
\end{equation}
Clearly, the vectors in~$V^\geq$ are all integer vectors in the outcome space
which are weakly dominated by some outcome vector $\bigl(f_1(\ve u),f_2(\ve
u),\dots, f_k(\ve u)\bigr)$ for a feasible solution $\ve x$ in $P\cap\Z^n$; however, we
have truncated away all outcome vectors which weakly dominate the
computed bound~$\ve{\bar v}$.  Let us consider the generating
function of~$V^\geq$, the multivariate polynomial
\begin{displaymath}
  g(V^\geq; \ve z)\
  = \sum \bigl\{\,  \ve z^{\ve v} : \ve v \in V^\geq \,\bigr\}.
\end{displaymath}
In the terminology of polynomial ideals, the monomials in $g(V^\geq;
\ve z)$ form a truncated ideal generated by the Pareto optima.  By the
Projection Theorem (our \autoref{project}), we can compute
$g(V^\geq;\ve z)$ in the form of a polynomial-size rational function
in polynomial time.\smallbreak

\noindent\emph{Step 2.} 
Let $\VPareto\subseteq\Z^k$ denote the set of Pareto optima.
Clearly we have
\begin{displaymath}
  \VPareto = \left(V^\geq \setminus (\ve e_1 + V^\geq) \right)
  \cap\dots\cap \left(V^\geq \setminus (\ve e_k + V^\geq) \right),
\end{displaymath}
where $\ve e_i\in\Z^k$ denotes the $i$-th unit vector and 
\begin{displaymath}
  \ve e_i + V^\geq = \{\, \ve e_i + \ve v : \ve v \in V^\geq \,\}.
\end{displaymath}
The generating function $g(\VPareto; \ve z)$ can be computed by the
Boolean Operations Lemma (\autoref{unioncomplement}) in polynomial time from $g(V^\geq; \ve z)$ as
\begin{equation}
  \begin{aligned}
    g(\VPareto; \ve z)={} &\bigl( g(V^\geq; \ve z) - g(V^\geq; \ve z) *
    z_1 g(V^\geq; \ve z) \bigr) \\
    &\quad * \dots * 
    \bigl( g(V^\geq; \ve z) - g(V^\geq; \ve z) *
    z_k g(V^\geq; \ve z) \bigr),
  \end{aligned}
\end{equation}
where $*$ denotes taking the Hadamard product of the rational functions.
\smallbreak

\noindent\emph{Step 3.} To obtain the number of Pareto optima, we compute the
specialization $g(\VPareto; \ve z\,{=}\,\ve 1)$.  This is possible in polynomial
time using residue techniques as outlined before the beginning of the proof.
\medbreak

\noindent\emph{Proof of part~(\ref{mainsec1-part-2}).} Now we recover the Pareto
\emph{strategies} that gave rise to the Pareto optima, i.e., we
compute a generating function for the set
\begin{displaymath}
  U^{\textrm{Pareto}} = \bigl\{\, \ve u \in\Z^n : 
  \text{$\ve u\in P\cap\Z^n$ and $\ve f(\ve u)$ is a Pareto optimum}\,\bigr\}.
\end{displaymath}
To this end, we first compute the generating function for the set
\begin{displaymath}
  S^{\textrm{Pareto}} = \bigl\{\, (\ve u, \ve v) \in\Z^n\times\Z^k: 
  \text{$\ve v$ is a Pareto point with Pareto strategy~$\ve u$}\,\bigr\}.
\end{displaymath}
For this purpose, we consider the multi-graph of the objective
functions $f_1,\dots,f_k$ over~$P$, 
\begin{equation}
  \label{eq:multi-graph}
  \begin{aligned}
    P^=_{f_1,\dots,f_k} = \bigl\{\, (\ve u, \ve v)\in\R^n\times\R^k :{}
    &\ve u\in P,\\
    & v_i = f_i(\ve u)\text{ for $i=1,\dots,k$} \,
    \bigr\}.
  \end{aligned}
\end{equation}
Using Barvinok's theorem, we can compute in polynomial time
the generating function for the integer points in~$P$, 
\begin{displaymath}
  g(P;\ve x)
  = \sum \bigl\{\,  \ve x^{\ve u} :
  \ve u \in P\cap \Z^n\,\bigr\},
\end{displaymath}
and also, using the monomial substitution $x_j
\rightarrow x_j z_1^{f_1(\ve e_j)} \cdots z_k^{f_k(\ve e_j)}$ for all~$j$,
the generating function is transformed into
\begin{displaymath}
  g(P^=_{f_1,\dots, f_k};\ve x, \ve z)
  = \sum \bigl\{\,  \ve x^{\ve u} \ve z^{\ve v} :
  (\ve u, \ve v) \in P^=_{f_1,\dots,f_k}\cap(\Z^n\times\Z^k)\,\bigr\},
\end{displaymath}
where the variables $\ve x$ carry on the monomial exponents the information of the
$\ve u$-coordinates of $P^=_{f_1,\dots,f_k}$ and the $\ve z$ variables of the
generating  function carry the $\ve v$-coordinates of lattice points in 
$P^=_{f_1,\dots,f_k}$. 
Now
\begin{equation}
  g(S^\textrm{Pareto};\ve x,\ve z) = 
  \bigl(g(P;\ve x)\, g(\VPareto;\ve z)\bigr)
  * g(P^=_{f_1,\dots,f_k}; \ve x, \ve z),
\end{equation}
which can be computed in polynomial time for fixed dimension by the
theorems outlined early on this section.  Finally, to obtain the generating
function $g(U^{\mathrm{Pareto}};\ve x)$ of the Pareto strategies, we need to compute the
projection of~$S^{\mathrm{Pareto}}$ into the space of the strategy variables~$\ve u$.
Since the projection is one-to-one, it suffices to compute the specialization 
\begin{displaymath}
  g(U^{\mathrm{Pareto}};\ve x) = g(S^\textrm{Pareto};\ve x,\ve z\,{=}\,\ve1),
\end{displaymath}
which can be done in polynomial time.
\end{proof}

\section{Efficiently listing all Pareto optima}
\label{s:listing}

The Pareto optimum that corresponds to the ``best'' compromise between the
individual objective functions is often chosen in an \emph{interactive mode},
where a visualization of the Pareto optima is presented to the user, who
then chooses a Pareto optimum.  Since the outcome space frequently is of a too
large dimension for visualization, an important task is to list (explicitly
enumerate) the elements of the \emph{projection} of the Pareto set into some
lower-dimensional linear space.

It is clear that the set of Pareto optima (and thus also any
projection) is of exponential size in general,
ruling out the existence of a polynomial-time enumeration algorithm.  
In order to analyze the running time of an enumeration algorithm, we must turn
to \emph{output-sensitive complexity analysis}.

Various notions of output-sensitive efficiency have appeared in the
literature; we follow the discussion
of~\cite{johnson-yannakakis-papadimitriou-1988}.  Let $W\subseteq\Z^p$~be a
finite set to be enumerated.  An enumeration algorithm is
said to run in \emph{polynomial total time} if its running time is bounded by
a polynomial in the encoding size of the input and the output.  A stronger
notion is that of \emph{incremental polynomial time}: Such an algorithm
receives a list of solutions $\ve w_1,\dots,\ve w_N\in W$ as an additional
input.  In polynomial time, it outputs one solution $\ve w\in W
\setminus\{\ve w_1,\dots,\ve w_N\}$ or asserts that there are no more
solutions.  An even stronger notion is that of a \emph{polynomial-delay}
algorithm, which takes only polynomial time (in the encoding size of the
input) before the first solution is output, between successive outputs of
solutions, and after the last solution is output to the termination of the
algorithm.  Since the algorithm could take exponential time to output all
solutions, it could also build exponential-size data structures in the course
of the enumeration.
This observation gives rise to an even stronger notion of efficiency, a
\emph{polynomial-space polynomial-delay} enumeration
algorithm.

We also wish to prescribe an \emph{order}, like the lexicographic order, in which the
elements are to be enumerated.   We consider term orders $\prec_R$
on monomials $\ve y^{\ve w}$ that are defined as
in~\cite{mora-robbiano-1988}
by a non-negative integral $p \times p$-matrix $R$ of full rank.  Two 
monomials satisfy $\ve y^{\ve w_1}\prec_R 
\ve y^{\ve w_2}$ if and only if $R\ve w_1$ is lexicographically smaller
than $R\ve w_2$.  In other words, if $\ve r_1,\ldots,\ve r_n$ denote the
rows of~$R$, there is some $j\in\{1,\ldots,n\}$ such that
$\langle \ve r_i, \ve w_1\rangle = \langle \ve r_i,\ve w_2\rangle$ for
$i<j$, and $\langle \ve r_j,\ve w_1\rangle < \langle\ve r_j,\ve w_2\rangle$. For 
example, the unit matrix $R=I_n$ describes the lexicographic term ordering.

We prove the existence of a polynomial-space polynomial-delay prescribed-order
enumeration algorithm in a general setting, where the set $W$ to be enumerated
is given as the projection of a set presented by a rational generating function.

\begin{theorem} \label{Theorem: Output sensitive enumeration}
Let the dimension $k$ and the maximum number~$\ell$ of binomials in the
denominator be fixed.

Let $V\subseteq\Z^{k}$ be a bounded set of lattice points with
$V\subseteq[-M,M]^{k}$, given only by the bound $M\in\Z_+$ and its rational
generating function encoding $g(V;\ve z)$ with at most~$\ell$ binomials in each
denominator.  Let 
\[
W = \{\, \ve w\in\Z^p: \exists \ve t\in\Z^{k-p} 
\text{ such that } (\ve t,\ve w)\in V \,\}
\]
denote the projection of $V$ onto the last $p$~components.  
Let $\prec_R$ be the term order on monomials in $y_1,\ldots,y_p$ induced
by a given matrix~$R\in\N^{p\times p}$.

There exists a polynomial-space polynomial-delay enumeration algorithm for
the points in the projection~$W$, which outputs the points of~$W$ in the
order given by~$\prec_R$.   The algorithm can be implemented without using the
Projection Lemma. 
\end{theorem}

We remark that \autoref{Theorem: Output sensitive enumeration} is a stronger
result than what can be obtained by the repeated application of the
monomial-extraction technique of Lemma~7 from~\cite{latte2}, which would only
give an incremental polynomial time enumeration algorithm.

\begin{proof}
  We give a simple recursive algorithm that is based on the iterative bisection of
  intervals. 
  \begin{quote}
    \emph{Input:}  Lower and upper bound vectors $\ve l, \ve
    u\in\Z^p$.\par
    \noindent\emph{Output:} All vectors~$\ve w$ in~$W$ with $\ve l \leq 
    R \ve w \leq \ve u$, sorted in the order $\preceq_R$.\smallskip\par
    \begin{enumerate}[1.]
    \item If the set $W\cap \{\, \ve w : \ve l \leq R \ve w \leq \ve u \,\}$
      is empty, do nothing.
    \item Otherwise, if $\ve l = \ve u$, compute the unique
      point~$\ve w\in\Z^k$
      with $R\ve w = \ve l = \ve u$ and output~$\ve w$.
    \item 
      Otherwise, let $j$ be the smallest index with $l_j \neq u_j$.  We bisect
      the integer interval $\{l_j, \dots, u_j\}$ evenly into $\{l_j, \dots,
      m_j\}$ and $\{m_j+1, \dots, u_j\}$, where $m_j = \floor{\frac {l_j +
          u_j}2}$.  We invoke the algorithm recursively on 
      the first part, then on the second part, using the corresponding lower
      and upper bound vectors. 
    \end{enumerate}
  \end{quote}
  We first need to compute appropriate lower and upper bound vectors $\ve
  l,\ve u$ to start the algorithm.  To this end, let $N$ be the largest number
  in the matrix~$R$ and let $\ve l = - p M N\ve 1$ and $\ve u = p M N\ve 1$.
  Then $\ve l \leq R \ve w \leq \ve u$ holds for all~$\ve w\in W$.
  Clearly the encoding length of $\ve l$ and $\ve u$ is bounded polynomially
  in the input data.
  
  In step~1 of the algorithm, to determine whether 
  \begin{equation}
    \label{eq:projection-empty}
    W\cap \{\, \ve w : \ve l \leq R \ve w \leq \ve
    u \,\} = \emptyset,
  \end{equation}
  we consider the polytope
  \begin{equation}
    Q_{\ve l,\ve u} = [-M, M]^{k-p} \times \{\,\ve w \in \R^p : \ve l \leq R \ve w \leq \ve
  u \,\} \subseteq \R^k, 
  \end{equation}
  a parallelelepiped in~$\R^k$.  Since $W$ is the projection of~$V$ and since
  $V\subseteq[-M,M]^k$, we have \eqref{eq:projection-empty} if and only if
  $V\cap Q_{\ve l,\ve u} = \emptyset$.  The rational generating
  function~$g(Q_{\ve l,\ve u}; \ve z)$ can be computed in polynomial time.  By
  using the Intersection Lemma, 
  we can compute the rational generating function~$g(V\cap Q_{\ve l,\ve u};
  \ve z)$ in polynomial time.  The specialization $g(V\cap Q_{\ve l,\ve u};
  \ve z=\ve1)$ can also be computed in polynomial time.  It gives the number
  of lattice points in~$V\cap Q_{\ve l,\ve u}$; in particular, we can decide
  whether $V\cap Q_{\ve l,\ve u}=\emptyset$.
  
  It is clear that the algorithm outputs the elements of~$W$ in the order
  given by $\prec_R$.  We next show that the algorithm is a polynomial-space
  polynomial-delay enumeration algorithm.   The subproblem in step~1 only
  depends on the input data as stated in the theorem
  and on the vectors~$\ve l$ and $\ve u$, whose encoding length only decreases
  in recursive invocations.  Therefore each of the subproblems can be solved
  in polynomial time (thus also in polynomial space). 
  
  The recursion of the algorithm corresponds to a binary tree whose nodes are
  labeled by the bound vectors $\ve l$~and~$\ve u$.  There are two types of
  leaves in the tree, one corresponding to the
  ``empty-box'' situation~\eqref{eq:projection-empty} in step~1, 
  and one corresponding to the ``solution-output'' situation in step~2.  Inner
  nodes of the tree correspond to the recursive invocation of the algorithm in
  step~3. 
  It is clear that the depth of
  the recursion is $\mathrm{O}(p \log (pMN))$, because the integer intervals are
  bisected evenly.  Thus the stack space of the algorithm is polynomially
  bounded.  Since the algorithm does not maintain any global data structures,
  the whole algorithm uses polynomial space only.

  Let $\ve w_i\in W$ be an arbitrary solution and let $\ve w_{i+1}$ be its
  direct successor in the order~$\prec_R$.  We shall show that the algorithm
  only spends polynomial time between the output of~$\ve w_i$ and the output
  of~$\ve w_{i+1}$.  The key property of the recursion tree of the algorithm
  is the following:
  \begin{equation}
    \text{\parbox{.8\linewidth}{%
        Every inner node is the root of a subtree that contains at least one
        solution-output leaf.}}
    \label{eq:all-branches-are-useful}
  \end{equation}
  The reason for that property is the test for
  situation~\eqref{eq:projection-empty} in step~1 of the algorithm.
  Therefore, the algorithm can visit only $\mathrm{O}(p \log (pMN))$ inner
  nodes and empty-box leaves
  between the solution-output leaves for $\ve w_i$ and $\ve w_{i+1}$.  For the
  same reason, also the time before the first solution is output and the time
  after the last solution is output are polynomially bounded.
\end{proof}

The following corollary, which is a stronger formulation of
\autoref{mainsec1}\,(\ref{mainsec1-part-3}), is immediate.

\begin{corollary}
  Let $n$ and $k$ be fixed integers.
  There exist polynomial-space polynomial-delay enumeration algorithms to 
  enumerate the set of Pareto optima of the multicriterion integer linear
  program~\eqref{eq:multicrit-ilp}, the set of Pareto strategies, or 
  arbitrary projections thereof in lexicographic order (or an arbitrary term order).
\end{corollary}

\begin{remark}
  We remark that \autoref{Theorem: Output sensitive enumeration} is of general
  interest.  For instance, it also implies the existence of a polynomial-space
  polynomial-delay prescribed-order enumeration algorithm for Hilbert bases of
  rational polyhedral cones in fixed dimension.  

  Indeed, fix the dimension~$d$ and let $C = \cone\{\ve b_1,\dots,\ve b_n\}
  \subseteq \R^d$ be a pointed rational polyhedral cone.  The \emph{Hilbert
    basis} of~$C$ is defined as the inclusion-minimal set $H\subseteq
  C\cap\Z^d$ which generates $C\cap\Z^d$ as a monoid.  For
  \emph{simplicial} cones~$C$ (where $\ve b_1,\dots,\ve b_n$ are linearly
  independent), \citet{barvinok-woods-2003} proved that one can compute the
  rational generating function $g(H;\ve z)$ (having a constant number of
  binomials in the denominators) of the Hilbert basis of~$C\cap\Z^d$ using the
  Projection Theorem.  The same technique works for non-simplicial pointed
  cones.    Now
  \autoref{Theorem: Output sensitive enumeration} gives a polynomial-space
  polynomial-delay prescribed-order enumeration algorithm.
\end{remark}

\section{Selecting a Pareto optimum using global criteria}
\label{s:selection}
  
Now that we know that all Pareto optima of a multicriteria integer linear
programs can be encoded in a rational generating function, and that they can be
listed efficiently on the output size, we can aim to apply selection
criteria stated by a user. 
The advantage of our setup is that when we optimize a global objective
function it guarantees to
return a Pareto optimum, because we evaluate the global criterion only on the
Pareto optima. Let us start with the simplest global criterion which generalizes
the use of the $\ell_1$ norm distance function:

\begin{theorem} \label{Theorem: Output sensitive enumeration polyhedral}
Let the dimension $k$ and the maximum number~$\ell$ of binomials in the
denominator be fixed.

Let $V\subseteq\Z^{k}$ be a bounded set of lattice points with
$V\subseteq[-M,M]^{n+k}$, given only by the bound $M\in\Z_+$ and its rational
generating function encoding $g(V;\ve z)$ with at most $\ell$ binomials in the
denominators.   

Let $Q\subseteq\R^k$ be a rational convex central-symmetric polytope with
$\ve0\in\interior Q$, given by its vertex or inequality description.
Let the polyhedral norm $\mathopen\|\cdot\mathclose\|_Q$ be defined using the Minkowski
functional 
\begin{equation}\label{eq:minkowski-functional}
  \| \ve y \|_Q = \inf \{\, \lambda \geq 0 : \ve y\in \lambda Q \,\}.
\end{equation}
Finally, let a prescribed point $\ve{\hat v}\in\Z^k$ be given.

\begin{enumerate}[\rm (i)]
\item There exists a polynomial-time algorithm to find a point $\ve v\in V$ that 
  minimizes the distance $d_Q(\ve v, \ve{\hat v}) = \| \ve v - \ve{\hat v}
  \|_Q$ from the prescribed point.\label{polyhedral-opt}
\item There exists a polynomial-space polynomial-delay enumeration algorithm
  for enumerating the points of~$V$ in the order of increasing distances
  $d_Q$ from the prescribed point~$\ve{\hat v}$, refined by an arbitrary term
  order~$\prec_R$ given by a matrix~$R\in\N^{k\times k}$.
  \label{polyhedral-order-enum}
\end{enumerate}
\end{theorem}

\autoref{th:intro-polyhedral-opt}, as stated in the introduction, is an
immediate corollary of this theorem.

\begin{proof}
  Since the dimension~$k$ is fixed, we can compute an inequality description 
  \begin{displaymath}
    Q = \{\, \ve y \in\R^k : A\ve y\leq \ve b\,\}
  \end{displaymath}
  of~$Q$ with $A\in\Z^{m\times k}$ and $\ve b\in\Z^k$ in polynomial time, if
  $Q$ is not already given by an inequality description.  Let $\ve v\in V$ be
  arbitrary; then 
  \begin{align*}
    d_Q(\ve{\hat v}, \ve v) &= \norm{\ve v - \ve{\hat v}}_Q \\
    &= \inf\bigl\{\, \lambda \geq 0 : \ve v - \ve{\hat v} \in \lambda Q \,\bigr\} \\
    &= \min\bigl\{\, \lambda \geq 0 : \lambda \ve b \geq A(\ve v - \ve{\hat v}) \,\bigr\}.
  \end{align*}
  Thus there exists an index
  $i\in\{1,\dots,m\}$ such that
  \begin{align*}
    d_Q(\ve{\hat v}, \ve v) &= \frac{(A\ve v)_i - (A\ve{\hat v})_i}{b_i};
  \end{align*}
  so $d_Q(\ve{\hat v}, \ve v)$ is
  an integer multiple of $1/b_i$.  Hence for every $\ve v\in V$, we have that
  \begin{equation}\label{eq:binary-search-granularity}
    d_Q(\ve{\hat v}, \ve v) \in \frac1{\lcm(b_1,\dots,b_m)}\Z_+,
  \end{equation}
  where $\lcm(b_1,\dots,b_m)$ clearly is a number of polynomial
  encoding size.  On the other hand, every $\ve v\in V$ certainly satisfies 
  \begin{equation}\label{eq:binary-search-largest-norm}
    d_Q(\ve{\hat v}, \ve v) \leq ka\bigl(M + \max\{\abs{\hat v_1}, \dots, \abs{\hat v_d}\}\bigr)
  \end{equation}
  where $a$ is the largest number in~$A$, which is also a bound of polynomial
  encoding size.

  Using Barvinok's algorithm, we can compute the rational generating function
  $g(\ve{\hat v} + \lambda Q; \ve z)$ for any rational $\lambda$ of polynomial
  enoding size in polynomial time.  We can also compute the rational
  generating function $g(V \cap (\ve{\hat v} + \lambda Q); \ve z)$ using the
  Intersection Lemma.  By computing the specialization $g(V \cap (\ve{\hat v}
  + \lambda Q); \ve z=\ve 1)$, we can compute the number of points in $V \cap (\ve{\hat v}
  + \lambda Q)$, thus we can decide whether this set is empty or not. 

  Hence we can employ binary search for the smallest $\lambda\geq 0$ such that
  $V \cap (\ve{\hat v} + \lambda Q)$ is nonempty.  Because
  of~\eqref{eq:binary-search-granularity} and
  \eqref{eq:binary-search-largest-norm}, it runs in polynomial time.  
  By using  the recursive bisection algorithm of \autoref{Theorem: Output
    sensitive enumeration}, it is then possible to construct one Pareto optimum
  in $V \cap (\ve{\hat v} + \lambda Q)$ for part~(\ref{polyhedral-opt}), or to construct a
  sequence of Pareto optima in the desired order for
  part~(\ref{polyhedral-order-enum}).
\end{proof}

\medbreak

Now we consider a global criterion using a distance
function corresponding to a non-polyhedral norm like the Euclidean
norm~${\mathopen\|\cdot\mathclose\|}_2$ (or any other $\ell_p$-norm for 
$1<p<\infty$). 
We are able
to prove a very strong type of 
approximation result, a so-called fully polynomial-time approximation scheme
(\FPTAS{}), in a somewhat more general setting.

\begin{definition}[\FPTAS]
  Consider the optimization problems 
  \begin{subequations}
    \begin{gather}
      \max\{\, f(\ve v) : \ve v \in V\, \}, \label{eq:generic-max-problem}\\
      \min\{\, f(\ve v) : \ve v \in V\, \}. \label{eq:generic-min-problem}
    \end{gather}
  \end{subequations}
  A \emph{fully polynomial-time approximation scheme (\FPTAS{})}
  for the maximization problem~$\eqref{eq:generic-max-problem}$ 
  or the minimization problem~$\eqref{eq:generic-min-problem}$, respectively, is
  a family $\{\, 
  \mathcal A_\epsilon : \epsilon \in\Q,\, \epsilon>0\,\}$ of approximation algorithms
  $\mathcal A_\epsilon$, each of which returns an \emph{$\epsilon$-approximation},
  i.e., a~solution $\ve v_\epsilon\in V$ with  
  \begin{subequations}
    \label{eq:epsilon-approximation}
    \begin{gather}
      \label{eq:max-epsilon-approximation}
      f(\ve v_\epsilon) \geq (1 - \epsilon) f^*
      \quad\text{where}\quad 
      f^* = \max_{\ve v\in V} f(\ve v),\intertext{or, respectively,}
      \label{eq:min-epsilon-approximation}
      f(\ve v_\epsilon) \leq (1 + \epsilon) f^*
      \quad\text{where}\quad 
      f^* = \min_{\ve v\in V} f(\ve v),
    \end{gather}
  \end{subequations}
  such that the algorithms $\mathcal A_\epsilon$ run in time polynomial in the
  input size and~$\frac1\epsilon$.
\end{definition}
\begin{remark}
  An \FPTAS{} is based on the notion of
  $\epsilon$-approximation~\eqref{eq:epsilon-approximation}, which gives an
  approximation guarantee relative to the value~$f^*$ of an optimal solution.  
  It is clear that this notion is most useful for objective functions~$f$ that
  are non-negative on the feasible region~$V$.  Since the approximation quality of
  a solution changes when the objective function is changed by an additive
  constant, it is non-trivial to convert an \FPTAS{} for a maximization
  problem to an \FPTAS{} for a minimization problem. 
\end{remark}

\begin{figure}
  \centering
  \ifpdf
    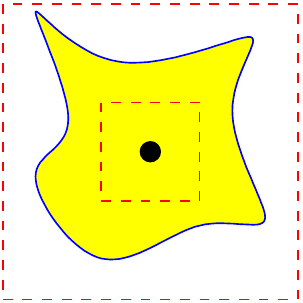
    \else
    \input{pseudonorm.pstex_t}
    \fi
  \caption{A set defining a pseudo-norm
    with the inscribed and circumscribed cubes $\alpha B_\infty$ and $\beta
    B_\infty$ (dashed).}
  \label{fig:pseudonorm}
\end{figure}

We shall present an \FPTAS{} for the problem of minimizing the distance of a
Pareto optimum from a prescribed outcome vector~$\ve{\hat v}\in\Z^k$.  We consider distances
$d(\ve{\hat v}, \cdot)$ induced by a pseudo-norm ${\mathopen\|\cdot\mathclose\|}_Q$ via
\begin{subequations}\label{eq:definition-pseudo-norm}
\begin{equation}
  d(\ve{\hat v}, \ve v) = {\| \ve v - \ve{\hat v} \|}_Q
\end{equation}
To this end, let $Q\subseteq\R^k$ be a compact basic semialgebraic set with 
$\ve 0\in\interior Q$, which is described by one polynomial inequality,
\begin{equation}
  Q = \bigl\{\, \ve y\in \R^k : q(\ve y) \leq 1 \,\bigr\},
\end{equation}
where $q\in\Q[y_1,\dots,y_k]$ is a homogeneous polynomial of (even) degree~$D$.  
The pseudo-norm ${\mathopen\|\cdot\mathclose\|}_Q$ is now defined using the
Minkowski functional 
\begin{equation}
  {\| \ve y \|}_Q = \inf \bigl\{\, \lambda \geq 0 : \ve y\in \lambda Q
  \,\bigr\} 
\end{equation}
\end{subequations}
Note that we do not make any assumptions of convexity 
of~$Q$, which
would make ${\mathopen\|\cdot\mathclose\|}_Q$ a~norm.  
Since $Q$ is compact and
$\ve0\in\interior Q$,
there exist positive rational numbers (norm equivalence constants) $\alpha$, $\beta$ with
\begin{equation}
  \label{eq:norm-equivalence}
  \alpha B_\infty \subseteq Q \subseteq \beta B_\infty
  \quad\text{where}\quad
  B_\infty = \bigl\{\, \ve y \in\R^k : {\|\ve y\|}_\infty \leq 1 \,\bigr\};
\end{equation}
see \autoref{fig:pseudonorm}.
\smallbreak

Now we can formulate our main theorem, which has
\autoref{th:intro-euclidean-fptas}, which we stated in the introduction, as an
immediate corollary. 
\begin{theorem}\label{th:pseudo-norm-fptas}
  Let the dimension~$n$ and the number~$k$ of objective functions be fixed.
  Moreover, let a degree~$D$ and two rational numbers $0<\alpha\leq\beta$ be fixed.
  Then there exists a fully polynomial-time approximation scheme for the problem of
  minimizing the distance $d_Q(\ve{\hat v}, \ve v)$, defined 
  via~\eqref{eq:definition-pseudo-norm} 
  by a homogeneous polynomial~$q\in\Q[y_1,\dots,y_k]$ of degree~$D$
  satisfying~\eqref{eq:norm-equivalence}, whose coefficients are encoded in
  binary and whose exponent vectors are encoded in unary,
  of a Pareto optimum of~\eqref{eq:multicrit-ilp} from a prescribed outcome vector
  $\ve{\hat v}\in\Z^k$.  
\end{theorem}

The proof is based on the following result, which is a more general
formulation of Theorem~1.1 from
\cite{deloera-hemmecke-koeppe-weismantel:intpoly-fixeddim}.
\begin{theorem}[\FPTAS{} for maximizing non-negative polynomials over finite
  lattice point sets] \label{th:ipo-fptas}
  For all fixed integers~$k$ (dimension) and~$s$ (maximum number of binomials
  in the denominator), there exists an algorithm with running 
  time polynomial in the encoding size of the problem and $\frac1\epsilon$ for
  the following problem.
  
  \textbf{Input:}
    Let $V\subseteq\Z^k$ be a finite set, given by a rational generating
    function in the form
    \begin{displaymath}
      g(V; \ve x) = 
      \sum_{i \in I} \gamma_i \frac { \ve x^{\ve c_i} }
      {  (1-\ve x^{\ve d_{i1}})  \dots (1-\ve x^{\ve d_{is_i}}) }
    \end{displaymath}
    where the the numbers $s_i$ of binomials in the denominators are at most~$s$.
    Furthermore, let two vectors $\ve v_{\mathrm L}$, $\ve v_{\mathrm
      U}\in\Z^k$ be given such that $V$ is contained in the box $\{\, \ve v:
    \ve v_{\mathrm L} \leq \ve v \leq \ve v_{\mathrm U}\,\}$.

    Let $f\in\Q[v_1,\dots,v_k]$ be a polynomial with rational coefficients
    that is non-negative on~$V$, given by a list of its
    monomials, whose coefficients are encoded in binary and whose exponents are
    encoded in unary. 

    Finally, let $\epsilon\in\Q$.

    \textbf{Output:}
    Compute a point $\ve v_\epsilon\in V$ that satisfies 
    \begin{displaymath}
      f(\ve v_\epsilon) \geq (1 - \epsilon) f^*
      \quad\text{where}\quad 
      f^* = \max_{\ve v\in V} f(\ve v).
    \end{displaymath}
\end{theorem}
In \cite{deloera-hemmecke-koeppe-weismantel:intpoly-fixeddim} the
above result was stated and proved only for sets~$V$ that consist of
the lattice points of a rational polytope; however, the same proof
yields the result above.

\begin{proof}[Proof of \autoref{th:pseudo-norm-fptas}]
  Using \autoref{mainsec1}, we first compute the rational generating
  function $g(\VPareto;\ve z)$ of the Pareto optima.
  With binary search using the Intersection Lemma with generating functions of
  cubes as in~\autoref{s:listing}, we can find the smallest non-negative
  integer~$\gamma$ such that 
  \begin{equation}
    \label{eq:smallest-nonempty-box}
    (\ve{\hat v} + \gamma B_\infty) \cap \VPareto \neq\emptyset.
  \end{equation}
  If $\gamma=0$, then the prescribed outcome vector $\ve{\hat v}$ itself is a
  Pareto 
  optimum, so it is the optimal solution to the problem.

  Otherwise, let $\ve v_0$ be an arbitrary outcome vector in $(\ve{\hat v} + \gamma B_\infty) \cap
  \VPareto$.  Then
  \begin{displaymath}
    \begin{aligned}
      \gamma \geq \norm{\ve v_0 - \ve{\hat v}}_\infty 
      &= \inf\set[big]{\,\lambda : \ve v_0 - \ve{\hat v} \in \lambda B_\infty \,} \\
      &\geq \inf\set[big]{\,\lambda : \ve v_0 - \ve{\hat v} \in \lambda
        \tfrac1\alpha Q \,} = \alpha \norm{\ve v_0 - \ve{\hat v}}_Q,
    \end{aligned}
  \end{displaymath}
  thus $ {\norm{\ve v_0 - \ve{\hat v}}}_Q \leq \gamma/\alpha$.
  Let $\delta = \beta\gamma/\alpha$.
  Then, for every $\ve v_1\in\R^k$ with $\norm{\ve v_1 - \ve{\hat
      v}}_\infty \geq \delta$ we have
  \begin{displaymath}
    \begin{aligned}
      \delta \leq \norm{\ve v_1 - \ve{\hat v}}_\infty 
      &= \inf\set[big]{\,\lambda : \ve v_1 - \ve{\hat v} \in \lambda B_\infty \,} \\
      &\leq \inf\set[big]{\,\lambda : \ve v_1 - \ve{\hat v} \in \lambda
        \tfrac1\beta Q \,} = \beta \norm{\ve v_1 - \ve{\hat v}}_Q,
    \end{aligned}
  \end{displaymath}
  thus 
  \begin{displaymath}
    \norm{\ve v_1 - \ve{\hat v}}_Q \geq \delta/\beta
    = \gamma/\alpha \geq \norm{\ve v_0 - \ve{\hat v}}_Q.
  \end{displaymath}
  Therefore, a Pareto optimum~$\ve v^*\in\VPareto$ minimizing the distance $d_Q$ from the
  prescribed outcome vector~$\ve{\hat v}$ is contained in the cube~$\ve{\hat v} +
  \delta B_\infty$.  Moreover, for all points $\ve v \in \ve{\hat v} + \delta
  B_\infty$ we have 
  \begin{displaymath}
    {\norm{\ve v_0 - \ve{\hat v}}}_Q \leq \delta/\alpha
    = \beta\gamma/\alpha^2.
  \end{displaymath}
  We define a function $f$ by 
  \begin{equation}
    f(\ve v) =  \paren{\beta\gamma/\alpha^2}^D - \norm{\ve v - \ve{\hat v}}_Q^D,
  \end{equation}
  which is non-negative over the cube $\ve{\hat v} + \delta
  B_\infty$. 
  Since $q$ is a homogeneous polynomial of degree~$D$, we obtain
  \begin{equation}
    f(\ve v) =  \paren{\beta\gamma/\alpha^2}^D - q(\ve v - \ve{\hat v})
  \end{equation}
  so $f$ is a polynomial.
  
  We next compute the rational generating function 
  \begin{displaymath}
    g(\VPareto\cap
    (\ve{\hat v} + \delta B_\infty);\ve z)
  \end{displaymath}
  from $g(\VPareto;\ve z)$ using the Intersection Lemma.  
  Let $\epsilon' > 0$ be a rational number, which we will determine later.
  By \autoref{th:ipo-fptas}, we compute 
  a solution $\ve v_{\epsilon'}\in\VPareto$ with
  \begin{equation*}
    f(\ve v_{\epsilon'}) \geq (1 - \epsilon') f(\ve v^*),
  \end{equation*}
  or, equivalently,
  \begin{equation*}
    f(\ve v^*) - f(\ve v_{\epsilon'}) \leq \epsilon' \, f(\ve v^*).
  \end{equation*}
  Thus,
  \begin{align*}
    [d_Q(\ve{\hat v}, \ve v_{\epsilon'})]^D - [d_Q(\ve{\hat v}, \ve v^*)]^D
    & =  \norm{\ve v_{\epsilon'} - \ve{\hat v}}_Q^D - \norm{\ve v^* - \ve{\hat v}}_Q^D\\
    & =  f(\ve v^*) - f(\ve v_{\epsilon'}) \\
    & \leq \epsilon' \, f(\ve v^*) \\
    & = \epsilon' \, \paren{
      \paren{\beta\gamma/\alpha^2}^D - \norm{\ve v^* - \ve{\hat v}}_Q^D }.
  \end{align*}
  Since $\gamma$ is the smallest integer with
  \eqref{eq:smallest-nonempty-box}
  and also $\norm{\ve v^* -
    \ve{\hat v}}_\infty$ is an integer, we have 
  \begin{displaymath}
    \gamma \leq \norm{\ve v^* -
      \ve{\hat v}}_\infty \leq \beta \norm{\ve v^* - \ve{\hat v}}_Q.
  \end{displaymath}
  Thus,
  \begin{align*}
    [d_Q(\ve{\hat v}, \ve v_{\epsilon'})]^D - [d_Q(\ve{\hat v}, \ve v^*)]^D
    & \leq \epsilon' \bracket{
      \paren{ \frac\beta\alpha }^{2D} -1 }
    \norm{ \ve v^* - \ve{\hat v}}_Q^D.
  \end{align*}
  An elementary calculation 
  yields
  \begin{align*}
    d_Q(\ve{\hat v}, \ve v_{\epsilon'}) - d_Q(\ve{\hat v}, \ve v^*)
    & \leq  \frac{\epsilon'}{D} \bracket{
      \paren{ \frac\beta\alpha }^{2D} -1 } d_Q(\ve{\hat v}, \ve v^*).
  \end{align*}
  Thus we can choose
  \begin{equation}
    \epsilon' = \epsilon\, D  \bracket{
      \paren{ \frac\beta\alpha }^{2D} -1 } ^{-1}
  \end{equation}
  to get the desired estimate.  Since $\alpha$, $\beta$ and~$D$ are fixed
  constants, we have $\epsilon' = \Theta(\epsilon)$.  Thus the computation of 
  $\ve v_{\epsilon'}\in\VPareto$ by \autoref{th:ipo-fptas} runs in time
  polynomial in the input encoding size and $\frac1\epsilon$.
\end{proof}

\begin{remark}
  It is straightforward to extend this result to also include the $\ell_p$
  norms for \emph{odd} integers~$p$, by solving the approximation problem
  separately for all of the $2^k=\mathrm{O}(1)$ shifted orthants $\ve{\hat v}
  + O_{\ve\sigma} = \{\, \ve v : \sigma_i(v_i-\hat v_i)\geq0\,\}$, where
  $\ve\sigma\in\{\pm1\}^k$.  On each of the orthants, the $\ell_p$-norm has a 
  representation by a polynomial as required by
  \autoref{th:pseudo-norm-fptas}. 
\end{remark}

\clearpage
\bibliographystyle{abbrvnat}
\bibliography{barvinok,iba-bib,weismantel}

\end{document}